\documentclass[12pt]{article}
\usepackage[all]{xy}
\UseComputerModernTips

\newtheorem{theorem}{Theorem}[section]
\newtheorem{proposition}[theorem]{Proposition}
\newtheorem{lemma}[theorem]{Lemma}

\newtheorem{definition}[theorem]{Definition}
\newtheorem{remark}[theorem]{Remark}
\newtheorem{example}[theorem]{Example}

\def\cc{{\cal C}}
\def\Z{{\bf Z}} 
\def\la{\longrightarrow}
\def\be{\begin{equation}}
\def\ee{\end{equation}}
\def\lotimes{\stackrel{L}{\otimes}}

\title{Braided $n$-categories and $\Sigma$-structures}
\author{Lawrence Breen\thanks{UMR CNRS 7539, Institut Galil\'ee, Universit\'e Paris 13, F-93430
Villetaneuse}}
\date{}

\begin{document}

\maketitle

\begin{abstract}
We associate to any braided 2-groupoid with vanishing intermediate homotopy group a principal bundle
(or torsor) endowed with a so-called
$\Sigma$-structure, and show that this is the natural generalization to the 2-category context of the
familiar quadratic invariant describing a braided groupoid. The corresponding structures for  higher
braided $n$-groupoids are also  examined. This leads  to the concept  of  $\Gamma_3$-torsor
pairs, which are in the same relation to cubic forms as torsors endowed with a $\Sigma$-structure are
to quadratic ones.  The discussion also covers the corresponding properties of braided $2$- and
$n$-stacks in groupoids.
\end{abstract}
\section{Introduction}
\setcounter{equation}{-1}%

\indent \indent This text is one of a series in which we explore the commutativity conditions which
may be imposed on tensor laws in monoidal $n$-categories. We showed  in
\cite{lb:moncat} that the structure of a monoidal groupoid
$\cc$ could be analyzed in terms of a simpler geometric object, the associated commutator 
biextension, which measures the obstruction to the full commutativity of the tensor law of $\cc$.
This commutator biextension, which mimicks the well-known commutator map associated to a central
extension of abelian groups, is
 endowed with a very strong anti-symmetry property, which we termed its alternating
structure \cite{lb:alt}. Such alternating biextensions are  classified up to isomorphism by the
groups
$\mathrm{Ext}^1(L\Lambda^2B,A)$, where $B$ and $A$ are a pair of abelian groups, and  $L\Lambda^2B$
is the (non-additive) derived functor of the second exterior power functor
$\Lambda^2_{\mathbf{Z}}$.  In \cite{lb:moncat}, we explored the manner in which monoidal
$n$-groupoids may be described  for $n \geq 2$ in terms of  higher analogs  of  alternating
biextensions, related  to the higher  derived exterior powers
$L\Lambda^i_{\mathbf{Z}}B$.  We always supposed there, as
we will here, that  the only nonvanishing homotopy groups of the
$n$-groupoids considered were the extremal ones, the group
$\pi_0(\cc)= B$ of isomorphism classes of objects and the group $\pi_n(\cc)= A$ of self-arrows of
the identity $(n$-$1)$-arrow, since without this hypothesis the associated geometric  objects
are no longer standard ones. Such
$n$-groupoids, together with their monoidal structure, are entirely determined up to equivalence by
the associated
$k$-invariant of their geometric realization, an element which lives in one of the cohomology groups
$H^{r+n}(K(B, r) , A)$. The integer $r$ is determined by the level of commutativity of the given
tensor multiplication law on$\cc$. 

\bigskip
 In the
present text, we consider monoidal
$n$-groupoids
$\cc$  endowed with a braiding structure, in other words  those   classified up to
equivalence by an element of the cohomology groups
$H^{n+2}(K(B,2),A)$. We  address  the  problem of determining the obstruction to the full
commutativity  of such a
monoidal structure.  In this situation, which is more restrictive than that examined in
\cite{lb:moncat} where no braiding axiom was postulated in $\cc$, the role previously assigned to
the alternating commutator map associated to a central  extension is taken up by the map
$\tau:
\Gamma_2^{\Z}B \la A$ introduced by A. Joyal and R. Street in 
\cite{js:braided}, whose source is the second divided power of the $\Z$-module $B$. The latter
associates to an object
$X$ of
$\cc$ whose class in 
$B$ is
$x$ the automorphism $\tau(x)$ of
$X
\otimes X$ determined by the braiding, viewed as
an element of  $A$. In the first situation examined here, that of a braided 2-groupoid $\cc$, this
 map
$\tau$ is replaced by the (trivial)  $A$-principal bundle (or torsor) on $B$ whose fiber above an
element
$x
\in B$ is the $A$-torsor of 2-arrows from the identity 1-arrow of $X \otimes X$ to the braiding
1-arrow $\tau(x)$. The additional structure which this $A$-torsor posesses asserts that it
satisfies the theorem of the cube, and that it is symmetric ({\it i.e.} invariant in a strong sense
under pull-back by the inverse map of $B$). Such $A$-torsors were said in \cite{lb:cube} to be
endowed with a $\Sigma$-structure. They are much more classical objects than the alternating
biextensions encountered above, since they correspond in an algebro-geometric context, when $A$ is
the multiplicative group $G_m$,  to  symmetric polarizations of the group scheme $B$. 

\bigskip

 In the higher categorical situations, one is led to    geometric
objects associated  to the higher divided power functors
$\Gamma_i^{\bf Z}$. Since braided 3-categories do not yield any suprisingly new multilinear
structures, we mainly consider those structures  associated to  braided 4-categories. To these
4-categories correspond somewhat exotic geometric structures, which may be thought of as higher
analogs of the torsors with $\Sigma$-structures defined by braided 2-categories. While
$\Sigma$-structure $A$-torsors on $B$ were classified by the group
$\mathrm{Ext}^1(L\Gamma_2B,A)$ associated to the derived functor $L\Gamma_2B$ of $\Gamma_2$, these
new structures, which we have termed
$\Gamma_3$-$A$-torsor pairs on $B$, are classified by the group
$\mathrm{Ext}^1(L\Gamma_3B,A)$.  These $\Gamma_3$-torsor pairs  therefore stand in the same
relation to cubic forms as the $\Sigma$-structure torsors do to quadratic ones. We unravel their
definition here, but do not give a complete proof of the fact that the structures in question are
indeed classified by this $\mathrm{Ext}$ group. Instead, we limit ourselves to the very strong
plausibility argument given by proposition
\ref{pr:last}.

\bigskip

While we will speak freely in this text of $n$-categories for rather large integers $n$, ($n\leq
5$)  will not make explicit  the requisite definitions of such $n$-categorical structures. Even
when this is not explicitly stated, our
$n$-categories will always be assumed to be $n$-groupoids, in other words categories in which, for each
$i < n$,
$i$-arrows are invertible up to an $i+1$-arrow, and $n$-arrows are strictly invertible.  The progress
recently made by various authors in achieving a workable definition of an
$n$-category, which is reflected in other texts \cite{BaDo}, \cite{Str} of these Proceedings,
justifies in our opinion this choice of the  
$n$-category  terminology for what are in fact  specific
homotopy types determined by 2-stage Postnikov fibrations. A reader unwilling to endorse this point
of view should consider
 this text as providing a very explicit description of  a natural filtration on certain 
cohomology groups 
$H^{n+2}(K(B,2),A)$. Nor  will we  make explicit the corresponding definition of an $n$-stack, and simply
refer to
\cite{lb:2-gerbes} for an illustration of what we have in mind in the case $n=2$. An $n$-stacks in
groupoids with trivial non extremal homotopy sheaves corresponds to a simplicial 2-stage
Postnikov simplicial sheaf,  described up to  equivalence by the  hyper-cohomology
groups
$\mathbf{H}^{r+n}(K(B,r),A)$. Here the  extremal homotopy sheaves $\pi_i(\cc)$ associated to $\cc$
for are defined  in the standard manner by sheafification  of  the naive homotopy presheaves
associated to
$\cc$. The  drawback with such a definition is that much of the
geometry underlying the idea of an $n$-category,  or  an $n$-stack, viewed as a family of objects
together with an  assortment of arrows recedes into the background.  We refer to
\cite{simpson},\cite{Hir-Sim} for an up to date discussion  of the theory of
$n$-stacks.  In the present context, the theory which we develop for $n$-groupoids extends to 
$n$-stacks in groupoids (which we will henceforth simply refer to as $n$-stacks). The  biextensions,
$\Sigma$-structures,  and other geometric structures which we encounter in the $n$-groupoid case
still exist, but   no longer posess global sections which trivialize the underlying torsors.
For the rest the situation is unchanged, with the proviso that additional non trivial terms in the
filtration of the cohomology must be considered,  which correspond to non vanishing higher
$\mathrm{Ext}$ groups.

\section{Braided categories  and stacks}
\setcounter{equation}{0}%
\indent \indent  We begin by reviewing the classification of braided groupoids.  Such a category may
be defined directly, as in \cite{js:braided}, as a groupoid endowed with a monoidal structure  for
which there exist inverses for objects and endowed with braiding morphisms\footnote{We henceforth
will denote by $XY$, rather than by the more customary $X
\otimes Y$, the product  of two objects $X$ and
$Y$ in a monoidal category 
$\cc$. On the other hand we will use the additive notation  for the composition law in $B$.}
\be 
\label{tauxy}
\tau_{X,Y}: YX \la XY \ee
 compatible with the
monoidal structure.  It may also be viewed as a 3-groupoid whose set of objets and of 1-arrows are
both reduced to a single element. The latter description makes it clear how to classify such
categories. Let 
$\cc$ be a braided groupoid.  The group
$\pi_0(\cc)$ of isomorphism classes of objects of such a braided groupoid $\cc$ and the group
$\pi_1(\cc)$
  of automorphisms
 of the identity object $I$ of $\cc$ are both abelian, and are respectively  denoted 
$B$ and $A$. The
monoidal structure on
$\cc$ determines for every object
$ X
\in
\cc$, by left multiplication by $X$, an identification 
\be
\label{timesX}
(X \otimes \, -) : B = \mathrm{Aut}_{\cc}(I)  \la \mathrm{Aut}_{\cc}(X)
\ee
 The nerve
$\mathcal{X}$ of
$\cc$ is a homotopy commutative $H$-space, which admits a double delooping to a space $\mathcal{Y} =
B^2\mathcal{X}$ whose only non-trivial homotopy groups are 
\[
\pi_i(\mathcal{Y})  = \left\{\begin{array}{ll} B & i=2 \\
A & i=3
\end{array} \right.
\] 
Such a space is entirely described, up to equivalence, by its unique $k$-invariant, an element
$k(\mathcal{Y})$ in the cohomology group $H^4(K(B, 2) , A)$. The entire $H$-space  structure of
$\mathcal{X}$, in other words the braided category structure of $\cc$, can be recovered from
$\mathcal{Y}$ by the double delooping $\mathcal{X}
\sim
\Omega^2\mathcal{Y}$, so that this invariant  yields a complete classification for equivalence
classes of  braided categories with associated groups $B$ and $A$. 

\medskip

Let us now examine the integral homology of $K(B,2)$ in low degrees. This computation, which is
functorial in $B$ and goes back to 
\cite{eml:hpin} II, yields the following result\be
\label{H(B,2)}
H_i(K(B,2))  = \left\{\begin{array}{ll} B & i=2 \\
0 & i=3 \\
\Gamma_2(B) & i = 4
\end{array} \right.
\ee
 A convenient method
for studying the integral homology of $K(B,2)$  in all degrees, along the lines of the computations
for
$K(B,1)$ in \cite{lb:fh}, consists in observing that the  homology of $K(B,2)$ (with the algebra
structure determined by the Pontrjagin product) is endowed with a natural divided power structure on
 even degree elements. The natural map
$B
\simeq H_2(K(B,2))$ determined by the fundamental class therefore extends to a functorial algebra
homomorphism
\be
\label{gammaH2} \Gamma_{\ast}(B) \la H_{\ast}(K(B,2))
\ee
When the abelian group $B$ is torsion-free, the entire homology of $K(B,2)$ is easy to compute.
The following result, which is implicit in \cite{cartan} (expos\'e 11, theorems 4 and 6), is
 explicitly stated  in  \cite{gd: thesis} (proposition 3.4).
\begin{lemma}
Let $B$ be a torsion-free abelian group. The canonical map (\ref{gammaH2}) is an isomorphism.
\end{lemma}
{\bf Proof:} 
For such a $B$,  the homology
$H_{\ast}(K(B,1))$ of the group
$B$ is isomorphic to the exterior algebra $\Lambda^{\ast}(B)$ of the ${\mathbf Z}$-module $B$. A
bar-construction (or Eilenberg-Moore spectral sequence) argument implies directly that 
\[ H_{\ast}(K(B,2)) \simeq \mathrm{Tor}^{\, \Lambda^{\ast}B}(\mathbf{Z},\mathbf{Z}) \]
 It is an immediate consequence of the
existence of the Koszul complex that the latter Tor term is isomorphic, as required, to  $ 
\Gamma_{\ast}B$.

\bigskip

The methods of \cite{lb:fh} now allow us to compute the homology of $K(B,2)$. Let $B_{\ast}
\la B$ be a torsion-free resolution of $B$. This determines a simplicial resolution $K(B_{\ast},0)$
of the constant simplicial group $K(B,0)$, as well as   a  
 torsion-free simplicial resolution $K(B_{\ast},2)$ of the simplicial abelian group
$K(B,2)$. By the previous lemma, the $E_1$ term of
the  homology spectral sequence associated as in  \cite{lb:fh} \S 1 to  $K(B_{\ast},2)$ can be identified
with the simplicial abelian group
$\Gamma_{\ast}(B_{\ast},0)$, obtained by applying the functor $\Gamma_{\ast}$ componentwise to
$K(B_{\ast},0)$. Since the homology of $K(B_{\ast},2)$ coincides with that of $K(B,2)$, the following
proposition is obtained.
\begin{proposition}
\label{pr:degss}
There exists a functorial   spectral sequence
\[ E^2_{p,q} = L_p \Gamma_q(B,0) \Longrightarrow H_{p+2q}(K(B,2))\]
which degenerates at $E^2$.
\end{proposition}
The degeneracy of the spectral sequence was essentially proved in   \cite{lb:fh} (with a
different description of the initial term). The latter is here the $p$th left (non-additive) derived
functor of
$\Gamma_q$, as defined in \cite{dp}, applied to
$B$ placed in degree 0, which will henceforth be simply denoted $L_p\Gamma_qB$ rather than the morse
standard $L_p\Gamma_q(B,0)$. In low degrees, this proposition immediately yields the computations
(\ref{H(B,2)}), as well as the isomorphism 
\be
\label{H5B2}
H_5(K(B,2)) \simeq L_1\Gamma_2B
\ee
We refer to \cite{lb:cube} \S 9 for a discussion of the functor $L_1\Gamma_2B$, and its
identification with  the functor $R(B)$ which by \cite{eml:hpin} theorem 22.1 describes the group
$H_5(K(B,2))$.

\medskip

We now consider the cohomology of $K(B,2)$. It follows from the universal coefficient
theorem and the computation (\ref{H(B,2)}) that  
\be
\label{H4G}
H^4(K(B, 2) , A) \simeq
\mathrm{Hom}(\Gamma_2(B), A)
\ee
so that  braided categories $\cc$  with associated groups $B$ and $A$ are  classified up to
equivalence by  quadratic maps $\tau_{\cc}: B \la A$. As observed by  A. Joyal and R. Street
(\cite{js:braided} theorem 3.3), this identification, which follows from the explicit description of
cycles for
$H_{4}(K(B,n))$ (as in
\cite{eml:hpin} theorem 26.1), associates to a category $\cc$ with a braiding \ref{tauxy}
the map 
\be
\label{tau}
\tau = \tau_{\cc} : B \la A
\ee which sends an element $x \in B$  to the automorphism
$\tau_{{X_x}, {X_x}}$ of $X_x X_x$, for some representative object $X_x \in \cc $ of $x$. The fact
that $\tau_{\cc}$ is independent of this choice of representative, and that it is  quadratic,
follows directly from (\ref{H4G}). It is, however,  of some interest to verify directly the
quadraticity of this map in geometric terms. To do this conveniently, we begin by assuming that the
multiplication law in our braided category $\cc$ is strictly associative, a situation which 
can always be achieved by replacing $\cc$ by an equivalent strictly associative monoidal category
$\mathcal{D}$, to which the braiding may be transported. In this new category, the  commutative
hexagons which determine the braided structure reduce to triangles such as the triangles
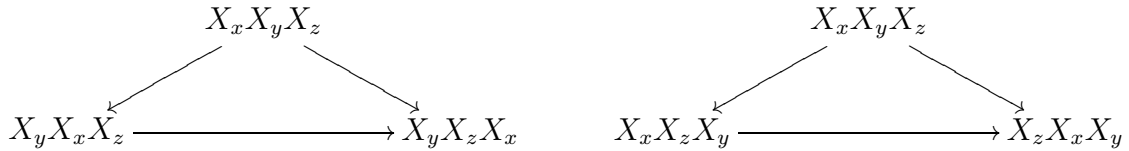
\[
\begin{array}{cccc}
{\xymatrix{&X_xX_yX_z \ar[dr] \ar[dl] & \\X_y X_xX_z \ar[rr] && X_y X_z X_x
}} & &&{\xymatrix{
& X_xX_yX_z \ar[dr] \ar[dl] & \\
X_x X_z X_y \ar[rr] && X_z X_x X_y}
}
\end{array}
\]
which respectively shuffle $X_x$ from  the left through $X_y X_z$ and  $X_z$ from the right
through $X_x X_y$. We denote the first of these `hexagon' diagrams by $H_{[x \mid_{{}_2}\,  y,z]}$
and the second by $H_{[x,y \, \mid_{{}_2} z]}$ and observe that these notations are consistent with
Eilenberg-Mac Lane's  for cells in $K(B,2)$. Let us choose, for each pair of elements
$x,y \in B$, an arrow
\be
\label{cxy}
\xymatrix{X_x X_y \ar[r]^{c_{x,y}} & X_{x+y}}
\ee
and henceforth denote, for any pair of elements $x,y \in B$, by $\tau_{x,y}$ the braiding
map $\tau_{{}_{X_x, X_y}}$ (\ref{tauxy}). We now consider  the commutative diagram
\be
\label{taudeg2}
\xymatrix{X_{x+y} \,  X_{x+y} \ar[r] \ar[dd]_{\tau (x+y)} & X_{x+y} \,  X_x X_y  \ar[dd] \ar[r]
\ar[dr] & X_x X_y X_x X_y \ar[dr]
\ar[r]^{\tau_{x,y}} & X_xX_xX_yX_y \ar[d]^{\tau (x)}\\
&&X_xX_{x+y} \ar[r] \ar[dl] & X_xX_x X_yX_y  \ar[dl] \ar[d]^{\tau (y)}Ê\\
X_{x+y} \, X_{x+y} \ar[r] & X_x X_y X_{x+y} \ar[r] & X_x X_yX_x X_y & X_x X_x X_y X_y
 \ar[l]_{\tau_{y,x}} }
\ee
Here the middle triangle is $H_{[x+y \mid_{{}_2}\,  x,y]}$, and the top and bottom right-hand ones
are respectively $H_{[x,y \, \mid_{{}_2} x]}$ and $H_{[x,y \, \mid_{{}_2} y]}$,  the
commutative squares being the obvious ones which follow from the functoriality of the maps
$\tau$. The commutativity of the outer square, when interpreted via arrows (\ref{cxy}) as a square
with the same object $X_{2x + 2y}$ at all vertices, asserts that the identity
\[
\tau(x+y) - \tau (x) -\tau (y) =  g(x,y) + g(y,x) \] is satisfied in $B$, where $g(x,y) \in B$ is
defined by the commutative diagram
\be
\label{gxy}
\xymatrix{
X_yX_x \ar[r]^{\tau_{x,y}} \ar[d]_{c_{y,x}} & X_x X_y \ar[d]^{c_{x,y}} \\
X_{x+y} \ar[r]_{g(x,y)} & X_{x+y}}
\ee
Asserting that the map $\tau $ (\ref{tau}) is of degree two is therefore equivalent to
showing that that the map 
$
\phi(x,y): B \times B \la B
$
defined by 
\be
\label{phi} \phi(x,y) = g(x,y) + g(y,x)
 \ee
 is bilinear. In order to prove this, consider first of all the  well-known ``Yang-Baxter''
hexagon\footnote{This hexagon (which is actually a hexagon when the associativity in the monoidal
category is no longer assumed to be strict) is also referred to in the physics litterature as the
``ABC = CBA'' identity.} associated to any three objects
$X,Y,Z$ in the braided category
$\mathcal{D}$. This can be built up as follows from a pair of triangles  $H_{[x \mid_{{}_2}\,  y,z]}$ and
$H_{[x \mid_{{}_2}\,  z,y]}$
\be
\label{YB1}
\xymatrix{
& XZY \ar[r] \ar@{-->}[drr] & ZXY \ar[dr] &\\
XYZ \ar[ur] \ar@{-->}[drr]
\ar[dr] 
&&&ZYX \\
&YXZ \ar[r] & YZX \ar[ur] &
}
\ee
 (the other method for constructing the identical Yang-Baxter
hexagon is by the scheme
\be
\label{YB2}
\xymatrix{
& XZY \ar[r]  & ZXY \ar[dr] &\\
XYZ \ar[ur] \ar@{-->}[urr]
\ar[dr] 
&&&ZYX \\
&YXZ \ar[r] \ar@{-->}[urr]& YZX \ar[ur] &
}
\ee
 the upper and lower commutative triangles now being $H_{[x,y \, \mid_{{}_2} z]}$ and $H_{[y,x \,
\mid_{{}_2} z]}$). In both instances, the commutativity of the central square follows from the
functoriality of the braiding map (\ref{tauxy})). Returning to the study of the map (\ref{phi}),
its linearity in the first variable may be obtained by considering the following  commutative
diagram, built up out from the  Yang-Baxter hexagon  associated to the three objects $X_x, X_y,
X_z$, to which have been glued a pair of triangles  $H_{[x,z \, \mid_{{}_2} y]}$,  $H_{[y
\mid_{{}_2}\,  z,x]}$ and a square whose commutativity follows, as always, from the functoriality
of the braiding map:
\be
\label{linphi}
\xymatrix{
&&X_xX_zX_y \ar[r]^{\tau_{z,x}} & X_zX_xX_y \ar[dr]^{\tau_{y,x}}&&\\
&X_xX_yX_z  \ar[dr]^{\tau_{y,x}}  \ar[ur]^{\tau_{z,y}}&&&X_zX_yX_x \ar[dr]^{\tau_{x,y}}& \\
X_xX_zX_y  \ar[ur]^{\tau_{y,z}} \ar[rr]_{\tau_{{}_{X_y,X_xX_z}}} && X_y X_x X_z \ar[r]^{\tau_{z,x}}
\ar[rrd]_{\tau_{{}_{X_xX_z, X_y}}} & X_y X_zX_x
\ar[ur]^{\tau_{z,y}} \ar[rr]_{\tau_{{}_{X_zX_x,X_y}}}&& X_zX_xX_y\\
&&&&X_xX_zX_y \ar[ur]_{\tau_{z,x}}&}
\ee
The arrows (\ref{cxy}) allows us to transform each of the eight vertices of the outer
diagram into the object $X_{x+y+z}$. The left hand ({\it resp.} right-hand) composite
arrow then becomes our map
$\phi (z,y)$ ({\it resp.} $\phi (x,y)$), and the lower composite one is $\phi (xz,y)$. The two
final  arrows $\tau_{z,x}$ encountered, when going around this diagram are both transformed into 
the maps
$g_{x,z}
$ (\ref{gxy}), with opposite orientations. Since the group $A$ of
automorphisms of an object $X$ in a braided category is abelian, these cancel out, and the
required identity 
\[\phi (xz,y) = \phi (x,y) + \phi (z , y) \]
is satisfied. The linearity of $\phi$ in the second variable $y$ follows automatically, since the
map $\phi$ is  symmetric by construction.

\bigskip

In order to verify the quadraticity of the map $\tau$ (\ref{tau}), it is still necessary to
ascertain that $\tau$  satisfies the identity 
\be
\label{tauquad}
\tau (x) = \tau (-x) \ee
Consider, for any $x \in B$, the commutative square 
\be
\label{tauquad2}
\xymatrix{
X_x X_x X_{-x}X_{-x} \ar[dr] \ar[r]^{\tau (-x)} \ar[d]_{\tau (x)} & X_x X_x X_{-x}X_{-x}
\ar[d]^{\tau_{x,-x}}\\ X_x X_x X_{-x}X_{-x} \ar[r]_{\tau_{x,-x}}& X_x X_{-x} X_x X_{-x}
}
\ee
built from an upper triangle $H_{[x,-x \, \mid_{{}_2} -x]}$ and a lower triangle $H_{[x
\mid_{{}_2}\,  x,-x]}$. The common diagonal, when  viewed  after contraction  to the
identity object $I$ of the middle term $X_xX_{-x}$ at its source and one or the other 
of the expressions $ X_x X_{-x}$ at its target, is essentially the identity map from
$X_x X_x
$ to itself. The triangles may therefore be viewed as  having a common diagonal edge along
which they may be attached as shown. By   the process previously applied  to the diagram
(\ref{linphi}), the outer square in (\ref{tauquad2}) transforms to one
in which all four vertices equal to $I= X_0$. The identity
\[ \tau(x) + \phi(x,-x) = \tau (-x) + \phi (x,-x) \] may be read off from its outer edges, so that
the identity (\ref{tauquad}) is proved. 
\begin{remark}
{\rm {\it i}) In a category in which the multiplication law is not strictly associative, it is
possible to carry out the above argument directly, without passing to a strictly associative
equivalent category $\mathcal{D}$.  A direct verification  that $H_4(K(B,2)) \simeq \Gamma_2(B)$
makes it clear that the only significant pentagonal cell to be added in that case to diagram  
(\ref{taudeg2}) is an associativity pentagon $P_{x,y,x,y}$, inserted at the common vertex
$X_xX_yX_xX_y$ appearing in both the top and bottom position. A pentagon
$P_{x,-x,x,-x}$ must similarly be added at the  vertex $X_x X_{-x} X_x X_{-x}$ to diagram
(\ref{tauquad2})

\medskip
\hspace {2cm}{\it ii}) It also follows from (\ref{H(B,2)}), and the universal coefficient theorem,
that  
\begin{eqnarray}
\label{H32}
H^3(K(B,2) , A) & \simeq & \mathrm{Ext}^1(B,A) \\
H^2(K(B,2), A) & \simeq & \mathrm{Hom}(B,A)
\end{eqnarray}
 and these relations have geometric interpretations analogous to that given above for 
(\ref{H4G}). The first implies, as mentioned in \cite{js:braided} remark 3.3, that the 
equivalences $F: \cc \la
\mathcal{D}$ between $\cc$ and another braided  category $\mathcal{D}$ equivalent to
it, when viewed up to natural equivalence, form
  a principal homogenous space
under the group
$\mathrm{Ext}^1(B,A)$ of auto-equivalences of
$\cc$. The second relation asserts that the set of natural
transformations between the identity equivalence $1_{\cc}: \cc \la \cc$ and some other fixed
equivalence $u: \cc
\la \cc$ is either empty, or a principal homogeneous space under  the group $\mathrm{Hom}(B,A)$ of
natural transformations between $1_{\cc}$ and itself. We will not make explicit the geometrical
realization of these assertions, but observe that this yields, together with the interpretation
given above of (\ref{H4G}),  a complete description up to 2-equivalence  of the  (symmetric
monoidal) 2-category of braided  categories with associated abelian groups $B$ and $A$, analogous
to the one given in \cite{js:braided} for the underlying 1-category.  
 }
\end{remark}

\bigskip

The previous discussion may also be carried out in the more general context of  a braided stack
$\cc$.
$B$ and
$A$ are now abelian groups in the topos $T$ considered, with the proviso that a section $x$ of $B$
is only  locally represented by an object $X_x$ in $\cc$.  The isomorphisms  (\ref{H32}) are
essentially unchanged in this new situation, despite the fact that $H^i(K(B,2) , A)$  now denotes
the $A$-valued hypercohomology groups  $\mathbf{H}^i(K(B,2) , A)$ of the simplicial object $ K(B,2)$
of
$T$. The higher
$\mathrm{Ext}^i$ groups no longer vanish for $i > 1$ in categories of abelian sheaves, as they did in
the category of abelian groups, so that the isomorphism  (\ref{H4G}) must now be replaced by
the
 exact sequence 
\be
\label{seq1}
\mathrm{0} \la \mathrm{Ext}^2(B, A) \la H^4(K(B,2), A) \la \mathrm{Hom}(\Gamma_2(B), A) \la
\mathrm{Ext}^3(B, A)  \ee
which is a low degree consequence of the universal coefficient spectral sequence
\[\mathrm{Ext}^p(H_q(K(B, 2)), A) \Longrightarrow H^{p+q}(K(B,2), A)\]
and of computation (\ref{H(B,2)}). The map $H^4(K(B,2), A) \la \mathrm{Hom}(\Gamma_2(B), A) $ is
essentially the one studied above. While it is {\it a priori} only defined locally,  the fact that
the explicit definition  of  the map (\ref{tau}) was independent of the choice of 
representing object $X_x$ for a section $x \in B$ ensures that these local versions of
$\tau_{\cc}$ are compatible. They therefore glue in the stack context to a genuine, globally
defined, sheaf map
$\tau_{\cc}: B
\la A$. While
 the map
 \be
\label{H4Gam2}
H^4(K(B,2), A) \la \mathrm{Hom}(\Gamma_2(B), A)
\ee is no longer injective, its kernel
$ \mathrm{Ext}^2(B, A)$ once again has a natural interpretation  as the group of
equivalence classes of strictly symmetric monoidal categories with associated groups $B$ and $A$
\cite{pd:sga4}. The quadratic maps from $B$ to $A$ which live in the image of the map
(\ref{H4Gam2}) thus still have, in the stack situation, a natural interpretation as the group of
obstructions to the existence, on a braided monoidal category $\cc$ with associated groups $B$
and $A$, of a strictly symmetric structure. This may be directly understood by observing that,
 the vanishing of the map (\ref{tau}) implies that  $\tau_{X,X}$ 
(\ref{tauxy})     is the identity map  for all objects $X$ in $\cc$. It then follows from
diagram (\ref{taudeg2})   that for any pair of objects $X$ and $Y$ in $\cc$ the composite map
\[ 
\xymatrix{
YX \ar[r]^{\tau_{x,y}} & XY \ar[r]^{\tau_{y,x}} & YX
}
\]
is also  the identity map, so that $\cc$ is indeed, as asserted,  a strictly symmetric monoidal
category.

\section{Braided 2-categories and $\Sigma$-structures}
\setcounter{equation}{0}%
\indent \indent Let us now pass   from the study of categories to that of 2-categories. Braiding
structures on such categories were introduced by Kapranov and Voevodsky \cite{kv:2-cat} (see also
 \cite{baeneu}, \cite{daystreet} definition 12, \cite{crans}).  We restrict ourselves here to
2-groupoids ({\it resp.} 2-stacks) whose intermediate homotopy group $\pi_1(\cc)$ is trivial, in
other words those for which, for any given pair of arrows $f,g: X \la Y$ between a pair of objects
$X, Y$, there exists  ({\it resp.} there exists locally) a 2-arrow $u: f \Longrightarrow g$
connecting
$f$ to $g$. The first of the remaining homotopy groups of such a braided 2-groupoid $\cc$ is  the
abelian group
$\pi_{0}(\cc) = B$ of isomorphism classes of objects of 
$\cc$. For a given object $X$ of $\cc$ (which we may choose to be the identity object $I$ of $\cc$),
the previous hypothesis ensures that the symmetric monoidal category
$\mathcal{A} = \mathrm{Aut}_{\cc}(X)$ is equivalent to the category with a unique object and an
abelian group $A$ of arrows. To make things specific, we may think of $A$ as the group of 2-arrows
\[u: 1_I \Longrightarrow 1_I \]
 from the identity 1-arrow $1_I: I \la I$  of the identity object $I$ to itself.  In the stack
case, both $B$ and $A$ are abelian sheaves, and   $\mathcal{A}$ is
equivalent as a  symmetric monoidal stack to the stack
$\mathrm{Tors}(A)$ of
$A$-torsors. As discussed in \cite{lb:2-gerbes} \S 8, the unique extra invariant necessary for a
full description of
$\cc$ as a braided 2-category or stack is ``Postnikov invariant'' of $\cc$, an
element $k$ in the cohomology (or hypercohomology group 
$H^4(K(B,2), \mathcal{A})$. Under our vanishing hypothesis for $\pi_1(\cc)$, this reduces to an
element of the more traditional cohomology (or hypercohomology) group $H^5(K(B,2), A)$. A direct 
analysis of this group {\it via} the   universal coefficient theorem is easy to achieve since the
values of the requisite homology groups of $B$ are know to us by (\ref{H(B,2)}) and  (\ref{H5B2}).
Let us consider instead, along the lines of the discussion carried out for monoidal 2-categories in
\cite{lb:moncat}, the filtration on the cohomology group $H^5(K(B,2), A)$ induced by the
filtration on homology determined by proposition \ref{pr:degss}. In the so-called punctual case in which
2-categories (rather than to 2-stacks) are considered, the only term in the filtration which contributes
something non trivial is the one involving
$\Gamma_2$, so that one obtains an isomorphism
\be
\label{H5G2}
 H^5(K(B,2), A) \simeq \mathrm{Ext}^1(L\Gamma_2B, A) \ee
Both left and right hand term can be subjected to a further anaysis via universal coefficient
theorem. To the  standard decomposition
\[0 \la \mathrm{Ext}^1(H_4(K(B,2), A) \la  H^5(K(B,2), A) \la \mathrm{Hom}(H_5(K(B,2)), A)
\la \mathrm{0}\]
of the left-hand term of (\ref{H5G2}) corresponds the decomposition
\[0 \la \mathrm{Ext}^1(\Gamma_2(B), A) \la  \mathrm{Ext}^1(L\Gamma_2B, A) \la
\mathrm{Hom}(L_1\Gamma_2(B), A) \la \mathrm{0}\]
of the right-hand one, and this correspondence  is consistent with the computation
(\ref{H(B,2)}), 
(\ref{H5B2}) of the low degree homology of $K(B,2)$.

\medskip

The  term
$\mathrm{Ext}^1(L\Gamma_2B, A)$ which appeared in (\ref{H5G2}) may be interpreted geometrically as
the group of  isomorphism classes of
$A$-torsors on
$B$ endowed with a so-called $\Sigma$-structure. For the reader's convenience we rapidly review this
concept, and refer to \cite{lb:cube} for additional information.  It is a strengthening of the more
familiar notion of a cube structure on such an
$A$-torsor  $P$ on $B$. The latter is determined by a trivialization $s$ of the ``second difference'' 
$\Theta (P)$ of the given torsor $P$, which is the $A$-torsor on $B^3$ defined by
\be
\label{deftheta}
 \Theta(P) = m_{123}^{\ast} P \wedge  (\bigwedge_{ij}m_{ij}^{\ast}P)^{-1} \wedge
(\bigwedge_i p_i^{\ast}P)
\ee 
where $\wedge$ is the contracted product of $A$-torsors, and the maps $m_{123}, m_{ij}, p_i$ from
$B^3$ to
$B$ are respectively the iterated sum, the partial
$(i,j)$th sum and the projection on the $i$th factor. In addition, the trivialization $s$ is required to
satisfy a cocycle condition, for which  we refer to  \cite{lb:cube} \S 2.  An alternate  method for
describing a cube structure on the $A$-torsor $P$ is to introduce its ``first difference'' of
the $A$-torsor $P$. This is the $A$-torsor $\Lambda (P)$ on  $B^2$ defined by 
\be
\label{fisrtdP}
\Lambda (P) = 
 m^{\ast}(P) \wedge p_1^{\ast}(P)^{-1} \wedge p_2^{\ast}(P)^{-1} 
\ee
(where $m$ is the group law and $p_i$ are the projections from $B \times B$ to $B$). The 
$A$-torsor $P$ is endowed with a cube structure if and only if
 $\Lambda (P)$ is not simply an $A$-torsor above $B^2$ but actually a biextension\footnote{ We refer
to 
\cite{ag:sga7} for a definition  of this concept.} of 
$B
\times B$ by
$A$.

\medskip

A $\Sigma$-structure on  an $A$-torsor
$P$ endowed with a cube structure is determined by the additional requirement that the torsor $P$ be
symmetric, in other words that there exists an isomorphism of $A$-torsors on $B$ 
\be
\label{sigmastr} \lambda : P \la i^{\ast}P 
\ee
(where $i: B \la B$ is the inverse map for the group law of $B$)  compatible in an appropriate sense
with the cube structure of
$P$ (see 
\cite{lb:cube}
\S 5).

\begin{example} 
{\rm In an
algebro-geometric context, when $A$ is the multiplicative group
$G_m$ and $B$ a commutative group variety, an $A$-torsor on $B$ is nothing else than a line bundle on
$B$. In particular, when the group $B$ is an abelian variety, we are assured of the existence of a
trivialization $s$ of
the induced line bundle $\Theta (L)$ by the classical theorem of the cube \cite{mu}. In
that situation,
 the cocycle
condition mentioned above is  satisfied, so that any such line bundle is automatically
endowed with a cube structure. Furthermore,   a symmetrization isomorphism
$\lambda$ (\ref{sigmastr}) is then always compatible  with the cube
structure, whenever $\lambda$ 
 respects the rigidification, in other words sends the point $s(e,e,e)$ in the fibre $P_e$ of $P$
above the unit element $e \in B$ to itself.}
   \end{example}

\bigskip Before explaining the correspondence (\ref{H5G2}) between braided 2-groupoids and
$\Sigma$-structures in more geometric terms, let us examine in purely cohomological terms the
corresponding situation in the 2-stack case. The identification (\ref{H5G2}) is now replaced by an
exact sequence 
\be
\label{seq:ex3}
 \rightarrow \mathrm{Hom} (\Gamma_2(B), A) \rightarrow \mathrm{Ext}^3(B, A) \rightarrow H^5(K(B,2),
A)
\stackrel{\tau}{\la }\mathrm{Ext}^1(L\Gamma_2B, A) \rightarrow \mathrm{Ext}^4(B, A)
\ee
which prolongs the exact sequence (\ref{seq1}). The term $\mathrm{Ext}^1(L\Gamma_2B, A) $ once more
describes $A$-torsors endowed with $\Sigma$-structures above $B$, but in the richer,
sheaf-theoretic context in which a torsor does not, as in the punctual case, always posess a global
section. The exactness of this sequence  shows that the  $\Sigma$-structure $ \tau
(\cc)$ associated by the middle map $\tau$ to a braided 2-stack $\cc$ is trivial if and only if
$\cc$ is Picard strict, in the terminology of \cite{lb:2-gerbes} definition 8.5.

\medskip

 This
Picard strict definition is the strongest possible commutativity which can be imposed on a 2-stack,
and it is of little interest in the punctual case  since the group $\mathrm{Ext}^3(B, A)$  which
 classifies such structures for abstract
abelian groups $B$ and $A$  vanishes. Let us briefly review from
\cite{lb:2-gerbes}, for the reader's convenience, the relations between this and other possible
commutativity conditions on 2-categories. This Picard strict commutativity condition is strictly
stronger than the condition on monoidal 2-categories which we called   ``strongly involutory'',
 which produces the 2-categories which R. Day and R. Street refer to as ``symmetric Gray
monoids''
\cite{daystreet}, and S. Crans as ``symmetric monoidal 2-categories''
\cite{crans}. The latter terminology for such 2-categories, which are classified under the vanishing
hypothesis  for
$\pi_1(\cc)$  by elements in the group
 $H^7(K(B,4),  A)$, seems the most appropriate one since the cohomology group in question lies in
the stable range. This  condition is  in turn   more restrictive than the condition 
which we refer to as ``weakly involutory'', and which applies to those categories which Day and
Street call ``sylleptic Gray monoids'' and  Crans refers to as ``sylleptic monoidal
2-categories''. The latter are classified, when $\pi_1(\cc)$ vanishes,  by elements in
$H^6(K(B,3),A)$. Finally, the commutativity conditions defining the original braided monoidal
2-categories of
\cite{kv:2-cat} are even weaker, since the latter are classified by elements of the
group $H^5(K(B,2), A)$. In each of these instances the universal coefficient theorem yields, in the
stack case, an edge-homomorphism map
\[\mathrm{Ext}^3(B, A) \la H^{i+3}(K(B,i), A)\]
and these maps are compatible with the suspension maps 
\[ H^{i+3}(K(B,i), A) \stackrel{S}{\la}   H^{i+2}(K(B,i-1), A) \] as $i$ varies.

\bigskip

If we now reason in geometric rather than homological terms, the $\Sigma$-structure associated
to a braided 2-stack $\cc$ may be understood by introducing the $A$-torsor  $\mathcal{S_C}$ on $B$,
whose sections above an section $x \in B$ are the 2-arrows 

\be
\label{SX}
\xymatrix{
X_xX_x \ar@/^1pc/[rr]^{1_{X_xX_x}} \ar@/_1pc/[rr]_{\tau(x)}& \: \: \Downarrow \:
\stackrel{S_x}{}  & X_xX_x }
\ee

 The $A$-torsor structure on
$\mathcal{S_C}$ is defined by composing such a locally defined 2-arrow with a 2-arrow from the
identity map
$1_{X_xX_x} $ to itself. In order to verify that  $\mathcal{S_C}$ is endowed with a
$\Sigma$-structure, we must first verify that it is endowed with a cube structure, in other words
that its ``first difference''
\[ \Lambda (\mathcal{S_C}) =  \mathrm{Isom} (p_1^{\ast}\mathcal{S_C} \wedge
p_2^{\ast}\mathcal{S_C}, m^{\ast}\mathcal{S_C}) \]
is a biextension of 
$B
\times B$ by
$A$. By definition, a local section of $\Lambda (\mathcal{S_C})$  is a rule which assigns to every
$(x,y)$ 
in some open set of $B \times B$, a triple  of 2-arrows  ($S_x, S_y, S_{x+y}$) (\ref{SX}), in a
manner consistent with the indicated action of $A$. Reinterpreting diagram  (\ref{taudeg2}) in a
2-categorical context, we observe that the commutative triangles  $H_{[x+y \mid_{{}_2}\,  x,y]}$, 
$H_{[x,y \,
\mid_{{}_2} x]}$ and
$H_{[x,y \, \mid_{{}_2} y]}$ have now been replaced by the corresponding hexagon 2-arrows
determined  by the braiding structure (see \cite {kv:2-cat}). The other more trivial
commutative cells in  diagram  (\ref{taudeg2}) are also replaced by  corresponding 2-arrows, so
that diagram (\ref{taudeg2}) determines, by composition of all these 2-arrows, a composite
2-arrow\footnote{This 2-arrow is roughly half of the diagram of 2-arrows   denoted
$((\bullet
\otimes
\bullet)
\otimes (\bullet \otimes \bullet))$  in \cite{kv:2-cat}, when specialized from
diagrams with vertices of the general form $XYZW$ and appropriate permutations thereof to those
for which $X=Z$ and $Y=W$.}
\be
\label{taudeg2-1}
\xymatrix{X_{x+y} \,  X_{x+y} \ar[r] \ar[dd]_{\tau (x+y)} &
  X_x X_y X_x X_y 
\ar[r]^{\tau_{x,y}} & X_xX_xX_yX_y  \ar@{}[ddll]^(0.43){\;}="1"
\ar@{}[ddll]^(0.57){\;}="2"
\ar[d]^{\tau (x)}\\
&  &  X_xX_x X_yX_y   \ar[d]^{\tau (y)} \\
X_{x+y} \, X_{x+y} \ar[r]   & X_x X_yX_x X_y & X_x X_x X_y X_y
 \ar[l]_{\tau_{y,x}}
\ar@{=>}"1";"2" }
\ee
 Composing it with the given  2-arrows ($S_x, S_y, S_{x+y}$)  according to the
scheme indicated in the diagram
\be
\label{taudeg2-2}
\xymatrix{X_{x+y} \,  X_{x+y} \ar[r] \ar@/_2pc/[dd]_1 \ar[dd]_{\Rightarrow}^{\tau (x+y)} &
  X_x X_y X_x X_y 
\ar[r]^{\tau_{x,y}} & X_xX_xX_yX_y 
\ar@{}[ddll]^(0.43){\;}="1"
\ar@{}[ddll]^(0.57){\;}="2"
\ar[d]_{\tau (x)}^{\Leftarrow} \ar@/^2pc/[d]^1\\
& &  X_xX_x X_yX_y   \ar[d]_{\tau (y)}^{\Leftarrow} \ar@/^2pc/[d]^1\\
X_{x+y} \, X_{x+y} \ar[r]   & X_x X_yX_x X_y & X_x X_x X_y X_y
 \ar[l]_{\tau_{y,x}}
\ar@{=>}"1";"2"  }
\ee
this yields a corresponding family of 2-arrows
\be
\label{SXY}
\xymatrix{
X_yX_x \ar@/^1pc/[rr]^{1_{X_yX_x}} \ar@/_1pc/[rr]_{\tau_{y,x} \, \tau_{x,y}}&\: \: \Downarrow \:
\stackrel{S_{x,y}}{}  & X_yX_x }
\ee

Let us now reconsider, in the 2-categorical context, the ``Yang-Baxter'' diagram (\ref{linphi}),
where all cells have now been replaced by 2-arrows determined by the braiding structure
of
$\cc$. It may be viewed as a rule which assigns to  every pair of 2-arrows 
(\ref{SXY}) 
\be
\label{SXYZ}
\begin{array}{ccc}
{\xymatrix{
X_xX_y \ar@/^1pc/[rr]^{1_{X_xX_y}} \ar@/_1pc/[rr]_{\tau_{x,y} \, \tau_{y,x}}&\: \: \Downarrow \:
\stackrel{S_{y,x}}{}  & X_xX_y }} & &{
\xymatrix{
X_zX_y \ar@/^1pc/[rr]^{1_{X_zX_y}} \ar@/_1pc/[rr]_{\tau_{z,y} \, \tau_{y,z}}& \: \: \Downarrow \:
\stackrel{S_{y,z}}{} & X_zX_y }
}
\end{array}
\ee
 a corresponding 2-arrow 
\be
\label{SXYZ1}
\xymatrix{
X_{xz}X_y \ar@/^1pc/[rr]^{1_{X_{xz}X_y}} \ar@/_1pc/[rr]_{\tau_{xz,y} \, \tau_{y,xz}}& \: \:
\Downarrow  \: \stackrel{S_{y,xz}}{}  & X_{xz}X_y }
\ee
 according the following scheme 
\be
\label{linphi2}
\xymatrix{
& \ar@{}[dd]^(.7){\;}="5"\ar@{}[dd]^(.9){\;}="6" &X_xX_zX_y \ar[r]^{\tau_{z,x}}_{\;}="8" & X_zX_xX_y
\ar[dr]_{\tau_{y,x}} 
\ar@/^2pc/[ddrr]^1_(.2){\;}="3" 
\ar@{}[ddrr]^(.3){\;}="4"
&\ar@{}[dd]^(.7){\;}="12"\ar@{}[dd]^(.9){\;}="13"&\\
 &X_xX_yX_z \ar[dr]_{\tau_{y,x}}  \ar[ur]_{\tau_{z,y}}&&&X_zX_yX_x \ar[dr]_{\tau_{x,y}}&
\\ X_xX_zX_y \ar@/^2pc/[uurr]^1_(.7){\;}="1" \ar@{}[uurr]^(.8){\;}="2"  \ar[ur]_{\tau_{y,z}}
\ar[rr]_{\tau_{{}_{X_y,X_xX_z}}} && X_y X_x X_z  \ar[r]^{\tau_{z,x}}="9"
\ar[rrd]_{\tau_{{}_{X_xX_z, X_y}}} & X_y X_zX_x
\ar[ur]_{\tau_{z,y}} \ar[rr]_{\tau_{{}_{X_zX_x,X_y}}}&& X_zX_xX_y\\
&&&&X_xX_zX_y \ar[ur]_{\tau_{z,x}}&
\ar@{=>}"1";"2"
\ar@{=>}"3";"4"
\ar@{=>}"5";"6"
\ar@{=>}"12";"13"
\ar@{}"8";"9"^(.4){\;}="10"
\ar@{}"8";"9"^(.6){\;}="11"
\ar@{=>}"10";"11"
}
\ee
Such a rule may be viewed as the second partial composition law in the biextension $\Lambda
(\mathcal{S_C})$, so that what remains to be verified is the associativity of this
law\footnote{According to  \cite{lb:cube} proposition 2.6, the partial group laws on a
symmetric biextension 
$\Lambda(L)$ associated to an $A$-torsor $L$ on $B$ endowed with a cube structures are
automatically  compatible. Only the commutativity and the associativity conditions need to be
verified.}. The commutativity is obtained by  reversing some arrows and  then reading diagram
(\ref{SXYZ1}) from right to left. The associativity asserts that for any given  triple of 2-arrows
$S_{y, w}, S_{y,z}, S_{y,x}, S_{y, w}$ constructed as above from given the corresponding set of
2-arrows
 (\ref{SX}), the two possible methods for constructing, by the method of 
(\ref{linphi2}), the corresponding 2-arrow $S_{y, wxz}$ 
\be
\label{SXYZW}
\xymatrix{
X_yX_{wxz} \ar@/^1pc/[rr]^{1_{X_yX_{wxz}}} \ar@/_1pc/[rr]_{\tau_{y,{wxz}} \, \tau_{wxz,y}}& \: \:
\Downarrow \: \stackrel{S_{wxz,y}}{}  & X_yX_{wxz} }
\ee
coincide. This  is a compatibility condition for the corresponding diagrams
(\ref{linphi2}) which we do not
display  here.

\bigskip

It remains to verify that the cube object $\mathcal{S_C}$ is actually endowed with a
$\Sigma$-structure. The requisite isomorphism (\ref{sigmastr}) between $\mathcal{S_C}$ and its
pullback
$i_B^{\ast}\,
\mathcal{S_C}$ by the inverse map $i_B$ of $B$ is easy to construct. It is  a rule which
associates  to a section 
$S_{-x}$ of the torsor $\mathcal{S_C}$ above $-x$   a section which lives above the opposite point
$x$ of $B$. This is determined by the pasting scheme

\be
\label{tauquad3}
\xymatrix{
X_x X_x X_{-x}X_{-x} \ar[dr] \ar[r]_{\tau (-x)}^{\:}="2" \ar@/^2pc/[r]^1_{\:}="1"
\ar[d]_{\tau (x)} & X_x X_x X_{-x}X_{-x}
\ar[d]^{\tau_{x,-x}} \ar@{}[dl]^(.2){\:}="3"  \ar@{}[dl]^(.4){\:}="4" 
\ar@{}[dl]^(.57){\:}="5" \ar@{}[dl]^(.77){\:}="6" \\ X_x X_x X_{-x}X_{-x}
\ar[r]_{\tau_{x,-x}}& X_x X_{-x} X_x X_{-x}
\ar@{=>}"1";"2"^{\; \;  \stackrel{S_{-x}}{}} 
\ar@{=>}"3";"4"
\ar@{=>}"5";"6"
}
\ee
 We do not attempt here a direct  verification of the 
compatibility between this rule and the ``Yang-Baxter'' rule, even though this is part of the
definition of a $\Sigma$-structure
\cite{lb:cube}.  

\begin{remark} \rm{It is now easy to understand in what sense $\mathcal{S_C}$ 
is an obstruction to commutativity of the group law of $\cc$: when  $\mathcal{S_C}$ is trivial, there
exists a global family of 2-arrows $S_x$. The axioms for a strict Picard 2-category, or
2-stack, as given in
\cite{lb:2-gerbes} definition 8.5, may now be read off as follows from the previous
$\Sigma$-structure conditions.  The requirement that  $S_{x}$ be additive with respect to $x$, in
the sense briefly described in \cite{lb:2-gerbes} page 150, corresponds in our context to the
assertion that the family of 2-arrows $S_{x,y}$ (\ref{SXY}) coincide with that determined by the
braiding structure on
$\cc$. The extra requirement that the 2-arrow displayed in
\cite{lb:2-gerbes} definition 8.5 coincide with $S_{x,x}$ is equivalent to the assertion above that
the composite 2-arrow  associated by (\ref{tauquad3}) to a given 2-arrow $S_{-x}$  is equal to 
$S_x$. The geometric content of the previous discussion may thus be summarized in very compact form
as the assertion that a braided 2-stack endowed with a family of 2-arrows $S_x$ satisfying the two
condition just mentioned, is a strict Picard 2-stack, {\it i.e.} automatically satisfies the
intermediate `Picard' commutativity axioms for a symmetric monoidal 2-stack. In  particular,
this asserts in the category context that a braided 2-category endowed with a family
of 2-arrows $S_x$ satisfying these conditions (and for which $\pi_1(\cc) = 0$) is
equivalent to the trivial symmetric monoidal 2-category  determined by a pair of abelian groups $B$
and
$A$. }
\end{remark}
\section {Comparison with the non-commutative obstructions}
\setcounter{equation}{0}%
\indent \indent It was shown in \cite{lb:moncat} \S 8 how to associate to any monoidal 2-category (or
2-stack)
$\cc$  determined as above by a pair of abelian groups $B$ and $A$ a so-called weak
(2,2)-extension $\mathcal{E}$ of $B \times B$ by $A$, which under an additional condition was found
to be alternating in a requisite sense. By the discussion at the very end of \cite{lb:moncat}, a
trivialization of  $\mathcal{E}$ compatible with this entire  structure makes $\cc$ into a
strict Picard stack (and therefore in the punctual case trivializes $\cc$ compatibly with its
monoidal structure).  We now compare this (2,2)-extension  obstruction $\mathcal{E}$ (which we 
will denote by $\mathcal{E}_{\cc}$ when we wish to emphasize its dependence on the  given 2-category
$\cc$) with the corresponding
$\Sigma$-structure obtruction $\mathcal{S}_{\cc}$ arising under the stronger
assumption that $\cc$ is braided.

 \bigskip

Before passing to this 2-category situation, it is   convenient to
examine the corresponding assertion for monoidal and braided 1-categories. As
recalled in \S 1, the obstruction to commutativity in a braided category or stack $\cc$  with
associated  invariant groups $B$ and $A$ is determined by the quadratic map $\tau_{\cc}$
(\ref{tau}) from
$B$ to
$A$. On the other hand, by  \cite{lb:moncat} proposition 3.1, there corresponds to any monoidal
category or stack $\cc$ with  associated invariants $B$ and $A$  an $A$-torsor
$ E_{\mathcal{C}} = E$ above $B \times B$, defined above $(x,y) \in B \times B$ by 
\be
\label{altE}
E_{x,y} = 
Isom_{\,\cc}(X_y X_x, X_x X_y) \ee
and which is endowed with the structure on a weak biextension\footnote{A weak biextension is simply
a (standard) biextension in which the pair of partial composition laws are no longer required to
be commutative, see \cite{ag:sga7} expos\'e VII 2.10.1.}. Furthermore, when $B$ is commutative, the
associated  commutator map, which examines whether the partial group laws on
$E$ are commutative, is a homomorphism $\phi_{\cc}$ from $\Lambda^3B$ to $A$, whose vanishing 
 ensures that $E$ is a (standard) biextension. The biextension $E$ is then automatically
an alternating one, a concept for which we refer to \cite{lb:moncat} \S 2,
\cite{lb:alt}. By
\cite{lb:moncat} corollaries 5.3 and 6.1, $E$ then measures the obstruction to the existence
on $\cc$ of a strict Picard structure. When the category $\cc$ is braided, the relation between the
quadratic map $\tau_{\cc}$ and  the alternating biextension $E_{\cc}$ goes as follows

\begin{proposition}
Let $\cc$ be a braided category with associated abelian groups $B$ and $A$.  The weak biextension
$ E_{\mathcal{C}}$ (\ref{altE}) associated to $\cc$ is a genuine biextension, and in fact a trivial
one. The quadratic map which by  \cite{lb:alt} proposition 1.5 determines the alternating structure
on this trivial biextension $E$ is, up to sign, the map $\tau_{\cc}$ (\ref{tau}).
\end{proposition}

{\bf Proof:} The braiding maps $\tau_{x,y}$ (\ref{tauxy}) determine a  section $s$ of the $A$-torsor
$E$ (\ref{altE}) above $B \times B$. The hexagon axioms assert that  this section is
bimultiplicative, and therefore trivializes $E$ as a biextension. This shows {\it a posteriori}
that  the weak biextension $ E_{\mathcal{C}}$ is in fact a genuine one. An alternating structure on
$E_{\cc}$ is determined by a trivialization of $\Delta E$ which differs from the trivialization
$\Delta s$ by a quadratic map $\tau$. The equation
\[ \tau(x) \circ \tau_{x,x}  = 1_{X_xX_x}\] in $\mathrm{Isom}_{\cc}(X_xX_x, X_x X_x)$ proves that
the the quadratic map $\tau$ associated to the canonical alternating structure  on $E$ determined as
in
\cite{lb:moncat} \S 5 by the identity arrow $1_{X_xX_x}$ on $X_xX_x$ is the opposite of the braiding
automorphism $\tau_{x,x}$ of  $1_{X_xX_x}$.

\begin{remark}  {\rm
This proposition translates in cohomological terms to the assertion that there exists, for any pair of abelian
groups  (or sheaves of abelian groups) $B$ and $A$, a diagram
\be
\label{diagG}
\xymatrix{
{\mathrm{0}} \ar[r] & \mathrm{Ext}^2(B, A) \ar[r] \ar@{=}[d] & H^4(K(B,2), A) \ar[r] \ar[d]_S &
{\mathrm{Hom}(\Gamma_2(B), A)}
\ar[d]_{\delta}
\\ {} & \mathrm{Ext}^2(B, A) \ar[r]  & G \ar[r]  &
\mathrm{Ext}^1(L\Lambda^2B, A) 
}
\ee
with exact lines, and for which the
left-hand square commutes, and the right-hand one anti-commutes. Here $G$ is the kernel of the
composite map
\be 
\label{ar:K}
H^3(K(B,1), A) \la \mathrm{Hom}(H_3(B) , A) \la \mathrm{Hom}(\Lambda^3B, A)
\ee
determined by the universal coefficient theorem and the iterated Pontrjiagin product map
$\Lambda^3B \la H_3(B)$. The vertical maps $S$ and $\delta$ are respectively the suspension map
and the boundary map induced on $\mathrm{Ext}^1$ by the distinguished triangle 
\be
\label{diag:gamlam}
L\Gamma_2B \la B \otimes B \la L\Lambda ^2B \stackrel{+1}{\la} L\Gamma_2B[1]
\ee
Finally the lower horizontal  exact sequence is described in detail, in geometric terms, in 
\cite{lb:moncat}. Its constituents are  low degree terms in the
 spectral sequence
\[E^{p,q}_2 = \mathrm{Ext}^p(L\Lambda^qB, A) \Longrightarrow H^{p+q}(K(B,1), A)\]
associated in cohomology to the homological spectral sequence 
(1.10) of \cite{lb:fh}.
}
\end{remark}

Passing now from categories to 2-categories, the situation is governed, when the intermediate
homotopy group $\pi_1(\cc)$ of the monoidal 2-category $\cc$ is trivial, by the diagram
 \be
\label{diagK}
\xymatrix{
{} \ar[r] & \mathrm{Ext}^3(B, A) \ar[r] \ar@{=}[d] & H^5(K(B,2), A) \ar[r] \ar[d]^S &
{\mathrm{Ext}^1(L\Gamma_2(B), A)}
\ar[d]^{\delta}
\\ {} \ar[r] & \mathrm{Ext}^3(B, A) \ar[r]  & K \ar[r]  &
\mathrm{Ext}^2(L\Lambda^2B, A) \\
}
\ee
Let us  briefly review the structure of the group
$H^4(K(B,1), A)$ which classifies such 2-categories $\cc$  \cite{lb:moncat} \S 8. This group is
endowed with a 4-term filtration, so that its subgroup $K$ displayed  above can no longer be
expressed (as the corresponding group $G$ in diagram (\ref{diagG})),   as the kernel of some 
edge-homomorphism analogous to the map (\ref{ar:K}). Instead, $K$ may be described as follows:
to  any  element of  
$H^4(K(B,1), A)$, representing a given 2-category $\cc$, is associated   a so-called weak 
(2,2)-extension $\mathcal{E_C} = \mathcal{E}$, whose fiber above a section $(x,y)
$ of  $B \times B$, is defined  by the formula 
\[
\label{2mock}
 \mathcal{E}_{x,y} = Isom_{\cc} (X_y X_x, X_x X_y)
\]
together with appropriate pair of partial tensor laws determined by the hexagon 2-arrows in $\cc$.
For the element in question to live in the subgroup $K$ of 
$H^4(K(B,1), A)$, it is necessary that a certain quadrilinear map $\psi_{\cc}: \Lambda^4B \la A$,
and  an induced alternating triextension $\Phi_{\cc} \in \mathrm{Ext}^1(L\Lambda^3B, A)$  both be
trivial. In that case 
$\mathcal{E}$ is endowed with  a genuine alternating (2,2)-extension structure, so that it
determines as required by diagram (\ref{diagK}) an element in
$\mathrm{Ext}^2(L\Lambda^2B, A)$.  As explained in
\cite{lb:moncat}
\S 8, the trivialization of this (2,2)-alternating biextension ensures that the monoidal 2-stack
$\cc$ is  strict Picard (and so is equivalent to the trivial monoidal category determined by the
pair of groups $B$ and $A$ in the 2-category case). Since such strict Picard stacks are classified
by elements of the group $\mathrm{Ext}^3(B,A)$, the present discussion gives a complete geometric
interpretation of the lower line of diagram (\ref{diagK}).  The upper line  in the diagram
was introduced in (\ref{seq:ex3}) above, where it was described in similar
geometrical terms. 

\medskip

 The comparison between the
$A$-torsor with $\Sigma$-structure $\mathcal{S_C}$ (\ref{SX})
 associated to a braided 2-stack $\cc$, and the alternating
(2,2)-biextension $\mathcal{E_C}$ associated to the underlying monoidal 2-stack  may now be carried
out along the same lines as in the corresponding discussion in the 1-category case. The given
braiding morphisms $\tau_{x,y}$ on $\cc$ determine a global section of the $A$-gerbe
$\mathcal{E_C}$. The  axioms satisfied by the hexagon 2-arrows in $\mathcal{C}$ imply that this
global section is compatible with the partial group laws on $\mathcal{E_C}$, and therefore
trivializes
$\mathcal{E_C}$ as a (2,2)-extension. In particular,  $\mathcal{E_C}$ is represented by a class which
lives in the subgroup $K$ of  $H^4(K(B,1), A)$, since this is where genuine (as opposed to weak)
(2,2)-biextensions live. The exactness of the  sequence 
\be 
\label{ex22}
 \cdots \rightarrow \mathrm{Ext}^1(L\Gamma_2(B), A) \rightarrow \mathrm{Ext}^2(L\Lambda^2B, A)
 \rightarrow \mathrm{Ext}^2(B \lotimes B, A) \rightarrow \mathrm{Ext}^2(L\Gamma_2(B), A) 
\rightarrow
\cdots
\ee
associated to the distinguished triangle  (\ref{diag:gamlam}) makes it clear that the alternating
structure on  $\mathcal{E_C}$ is determined by a trivialization of its restriction $\Delta
\mathcal{E_C}$ above the diagonal in $B$. Such a restriction to the diagonal of a (2,2)-extension
$\mathcal{E}$ is an
$A$-gerbe on
$B$ canonically endowed with a structure which we might call a 2-$\Sigma$ structure, and which is
classified up to equivalence by the group $\mathrm{Ext}^2(L\Gamma_2B, A)$.  The complete definition of such a 2-$\Sigma$ structure may
be obtained along the same lines as that of alternating (2,2)-extensions in \cite{lb:moncat}, by
examining the d\'evissage of this
$\mathrm{Ext}^2$ group associated to the distinguished triangle 
\[ L\mathrm{Sym}^2B \la L\Gamma^2B \la B \lotimes {\bf Z}/2 \stackrel{+1}{\longrightarrow } 
L\mathrm{Sym}^2B\: [1]\]
We do not work this out in detail here, since such a discussion was already carried out in a
related context in  \cite{lb:moncat} \S 8. We simply recall from
\cite{lb:moncat} that the  sought-after trivialization of 
$\Delta
\mathcal{E_C}$, which determines the alternating structure on $\mathcal{E_C}$, is given by the
identity map
\[1_{X_xX_x}: X_xX_x \la X_xX_x\] viewed as an object in $(\Delta
\mathcal{E_C})_x = \mathrm{Isom}(X_xX_x, X_xX_x)$ for a varying section $x \in B$.  The
$A$-torsor with $\Sigma$-structure which describes this alternating structure on the torsor 
$\mathcal{E_C}$ (trivialized by the given braiding arrows $\tau_{x,y}$ in 
$\cc$) is, up to a sign depending on our conventions,  the $A$-torsor $\mathcal{S_C}$ of (locally
defined) 2-arrows $S_x$ (\ref{SX}) which compare the  induced trivialization $\tau(x)$ of $\Delta
\mathcal{E_C}$ with the canonical trivialization $1_{X_xX_x}$.  It follows from this discussion
that the right-hand square in diagram (\ref{diagK}) commutes up to sign, where the vertical map
$\delta$ is defined either cohomologically as in (\ref{ex22}), or geometrically as above.

\section{Higher invariants: the  higher $\Sigma$-structure torsors}
\setcounter{equation}{0}%
\indent \indent  We now briefly examine the invariants associated to higher homotopy types. We think
of the homotopy types described by the groups
$H^6(K(B,2), A)$ as being represented by appropriately defined braided 3-categories $\cc$ (or
rather, 3-groupoids since we assume that every $n$-arrow in $\cc$ is invertible up to an
$n+1$-arrow) with the usual additional assumptions: the group $\pi_0(\cc)$ of isomorphism
classes of objects is isomorphic to $B$,
 the group
$\pi_3(\cc)$ of self 3-arrows  of any identity 2-arrow is isomorphic to $A$, and the  intermediate
homotopy groups
$\pi_i(\cc)$ are trivial. We will not spell out the braiding axioms here, but merely  observe 
that one method for obtaining them would be to think of $\cc$ as a 5-category with a single object
and a single 1-arrow. The multilinear objects attached, as in the lower dimensional case,  to
an appropriate filtration of $H^6(K(B,2), A)$ are now the groups $\mathrm{Hom}(\Gamma_3B,A),
\mathrm{Ext}^2(L\Gamma_2B, A)$ and
$\mathrm{Ext}^4(B,A)$. We have seen that the group $\mathrm{Ext}^2(L\Gamma_2B, A)$ describes classes
of what might be termed
$2$-$\Sigma$ structures, the higher (1-categorical) analog of the $\Sigma$-structure $A$-torsors
on 
$B$. We do not spell out here the definition of such structures, nor
of the strict Picard 3-categorical structures classified by  $\mathrm{Ext}^4(B,A)$ since this is
routine. Let us merely observe that 
$2$-$\Sigma$-structures are to (ordinary)
$\Sigma$-structures as their anti-symmetric analog, the alternating (2,2)-extensions of
\cite{lb:moncat} are to (ordinary) alternating biextensions. The only essentially new feature which
occurs when braided 3-categories are considered is thus the appearance of the degree three component
$\Gamma_3B$ of the divided power algebra $\Gamma_{\ast}B$ of $B$, together with the projection
\be
\label{UCT3}
H^6(K(B,2), A) \la \mathrm{Hom}(\Gamma_3B, A)
\ee
The latter map is readily understood from a topological point
of view, since it is induced by the isomorphism $\Gamma_3B \simeq H_6(K(B,2))$
 \cite{eml:hpin}\footnote{In \cite{eml:hpin}, the standard degree $i$ component
$\Gamma_iB$ of the divided power algebra $\Gamma_{\ast}B$ is somewhat confusingly denoted
$\Gamma_{2i}B$.}
theorem 21.1  determined by the divided power algebra structure of
$H_{\ast}(K(B,2)$, and by the universal coefficient theorem. In categorical terms, it can be
understood by observing that in a braided 3-category, the two possible methods of construction of a
Yang-Baxter 2-arrow described by diagrams (\ref{YB1},
\ref{YB2}) no longer yield the same result, as in a braided 2-category. Instead, they can only be
compared by a 3-arrow, which is part of the the braiding axioms for 3-categories, and which we
denote,  for 
  representative objects $X_x, X_y,X_z$ in
$\cc$ of the isomorphism classes $x,y,z
\in B$, by $(YB)_{x,y,z}$. Such a 3-arrow corresponds to the degree 6 cell in Eilenberg-Mac Lane's
chain complex for $K(B,2)$ traditionally denoted by $[\,x\! \!\mid_{{}_2} \! y \!\! \mid_{{}_2} 
\!z]$.
 It corresponds, in Getzler-Jones'  cell decomposition \cite{Getz-Jo} \S5.4 of the quotiented
configuration space
$\stackrel{\circ}{F_3}(2)$ of pairs of points in ${\mathbf R}^3$, to the cell
$[1\!\mid\mid\!2\!\mid\mid\!3]$. It now follows from the explicit computation of
the elements of $H_6(K(B,2))$ that the  map (\ref{UCT3})  associates to the given braided
3-category
$\cc$ the homorphism 
\be
\label{Ups}
\Upsilon_{\cc}: \Gamma_3B \la A
\ee
 which sends an element $\gamma_3(x)$ to the 3-arrow $(YB)_{x,x,x}$. The latter specialized
Yang-Baxter 3-arrow is, as required an element of the group $A$  of self 3-arrows in $\cc$ from
an  appropriate trivial 2-arrow to itself. 

\medskip

The group $\Gamma_3B$ differs from the previously encountered group $\Gamma_2B$, since it
is no longer generated by indecomposable elements of the form 
$\gamma_3(x)$,   additional generators of the form $\gamma_2(x)y$ with $x,y \in B$  being
required. The map (\ref{Ups}) can therefore no longer be entirely described by the values of the
elements
$(YB)_{x,x,x}
\in A$ for varying $x \in B$. Instead, one is also led to consider the  composite map 
\[ \Gamma_2B \otimes B \la \Gamma_3B \stackrel{\Upsilon_{\cc}}{\la} A \]
determined by the algebra structure in $\Gamma_{\ast}B$. The latter may  be viewed as a map  \[B
\times B \la A\]  which is quadratic in the first variable $x$ and linear in the second variable $y$.
It is constructed  by assembling the three Yang-Baxter
 3-arrows  $(YB)_{x,x,y}, (YB)_{x,y,x},(YB)_{y,x,x}$ into a self 3-arrow, in a manner which we will
not make explicit here.

\bigskip

We now pass to the groups $H^7(K(B,2), A)$, whose elements describes the  braided
4-categories with trivial intermediate homotopy groups. There are now three corresponding
multilinear invariants, which live in the groups
\[ \begin{array}{ccccc}
\mathrm{Ext}^1(L\Gamma_3B,A) & &
\mathrm{Ext}^3(L\Gamma_2B, A) && \mathrm{Ext}^5(B,A)
\end{array}
\]
 While the last two represent some
significant structure in the 4-stack case, they both vanish when 4-categories are considered, so we
will restrict ourselves to a discussion  of the geometric objects representing the remaining group
$\mathrm{Ext}^1(L\Gamma_3B,A)$. This is in any case the most interesting of the three, since it
involves the  degree three functor $\Gamma_3B$. This functor lives in the exact sequence 
\be
\label{syzg3} \Gamma_2B \otimes B \la \Gamma_3B \la B/3B \la \mathrm{0} \ee
reminiscent of the short exact sequence
\[\mathrm{0} \la \mathrm{Sym}^2B  \la \Gamma_2B \la B/2B \la \mathrm{0}\]
by which one analyzes a $\Sigma$-structures in terms of the associated symmetric biextension
\cite{lb:cube} \S 8. In the present situation, the first arrow is  the multiplication map 
considered  above, and the second the map  which sends $\gamma_3(x)$ to the class of $x  \: 
\mathrm{mod}\:(3B)$.

\medskip

Since there exist two sorts of generators for $\Gamma_3B$, the geometric structure classifed
by the group $\mathrm{Ext}^1(L\Gamma_3B,A)$ involves a pair of objects $(E, L)$ where $E$ an
$A$-torsor on $B \times B$  and $L$ is an $A$-torsor on $B$. When such a pair $(E, L)$ arises from a
braided 4-category $\mathcal{D}$ with trivial intermediate homotopy groups, the corresponding 
$A$-torsor $L_{\mathcal{D}}$
 is defined by 
\[L_{\mathcal{D}} = \mathrm{Isom}_{\mathcal{D}}(1 , \Delta^{\ast}_{123} (YB))\]
Its fiber above $x \in B$ is the torsor of 4-arrows from the identity 3-arrow to $(YB)_{x,x,x}$.
The fibre above ${x,y} \in B^2$ of the corresponding $A$-torsor $E_{\mathcal{D}}$ similarly consists
of the torsor of 4-arrows in $\mathcal{D}$ from the identity 3-arrow to the 3-arrow assembled from
$(YB)_{x,x,y}, (YB)_{x,y,x}$ and $(YB)_{y,x,x}$ to which we previously alluded (without making it
explicit) in the 3-category case.

\medskip

 Returning to the general situation, we now describe the additional structure to be imposed on
the  pairs $(E,L)$. The map from
$\mathrm{Ext}^1(\Gamma_3B,A)$ to
$\mathrm{Ext}^1(L\Gamma_2B \lotimes B, A)$ induced the multiplication in $\Gamma_{\ast}B$  is the
first projection from $(E,L)$ to $E$, so that the first element of structure to be imposed on the
pair
$(E,L)$ is that $E$ be quadratic in the first variable and linear in the second.  This means that
for each $y \in B$, $E_{(-,y)}$, viewed as a torsor above $B$, is endowed with a
$\Sigma$-structure, and that $E_{(x,-)}$, viewed for each $x \in B$ as an $A$-torsor on $B$, is a
commutative extension of  $B$ by $A$. We also require the compatibility of these two structures,
in other words that the partial multiplication map 
\be
\label{compat}
E_{x, 	y_1} \wedge E_{x, 	y_2} \la E_{x, 	y_1 +y_2}
\ee
determined by the second variable group law
respect the $\Sigma$-structures.

\medskip

The next set of conditions  on the pair $(E,L)$ arises from the relations
\begin{eqnarray}
\gamma_3(x+x') - \gamma_3(x) - \gamma_3(x') & = & \gamma_2(x) \, x' + \gamma_2(x') \,x \\
3 \gamma_3(x) & = & \gamma_2(x) \, x \label{gama3}
\end{eqnarray}
satisfied by the divided powers.
The first of these relations translates into the requirement that there exists an isomorphism of
$A$-torsors
\be
\label{alph}
\Lambda(L) \stackrel{\alpha}{\la} E \wedge s^{\ast}E
\ee
where $s$ is the map which permutes the factors in $B^2$, and $\Lambda(L)$ is  the
``first difference'' $A$-torsor of $L$, defined above $B^2$   formula (\ref{fisrtdP}). Similarly,
the relation (\ref{gama3}) translates to the requirement that there exists an isomorphism  of
$A$-torsors
\be
\label{bet}
\Delta E \stackrel{\beta}{\la}  L^3
\ee
between the restriction of $E$ to the diagonal and the contracted third power $L \wedge L \wedge L$ 
of the torsor $L$. It is obvious that such a requirement is satisfied by the pair $(E_{\mathcal{D}},
L_{\mathcal{D}})$ since sections of the restriction of $E_{\mathcal{D}}$ to the diagonal involve
three copies of $(YB)_{x,x,x}$.

\medskip

 Consider the pullback 
\[
\xymatrix{(\Lambda \times 1)^{\ast}(E \wedge s^{\ast} E) \ar[rr]^>>>>>>>>>>{{(\Lambda \times
1)}^{\ast} \alpha} && (\Lambda \times 1 )^{\ast} (\Lambda(L))}
\]
by $(\Lambda \times 1)$ of the map  $\alpha$ (\ref{alph}). The target of this map is the ``first
difference of the first difference '' of $L$, which is canonically isomorphic to  the ``second
difference''
$\Theta(L)$ of $L$ (\ref{deftheta}). By linearity of $E$ in the first variable, the source of
this map simplifies so that it reduces to an isomorphism
\be
\label{lambalph}
\xymatrix{
(\Lambda \times 1)^{\ast} E \ar[r] & {\Theta (L)}
}
\ee
Since $E$ posesses a $\Sigma$-structure with respect to its first variable, and is linear with
respect to the second one, the source of this map is
  canonically endowed with a triextension structure above $B \times B \times B$, and in fact one
which is symmetric with respect to the interchange of the first two variables. In particular the
$A$-torsor $L$ is endowed {\it via} this isomorphism with a hypercube structure, since its second
difference is multilinear\footnote{so that its third difference is trivial.}. 

\medskip

We now list  the conditions which the maps $\alpha$ and $\beta$ must satisfy. Observe that the
target $\Theta (L)$ of the map (\ref{lambalph}) is invariant up to canonical isomorphism under any
permutation of the three factors in
$B^3$. The canonical isomorphism $\Theta (L) \la \sigma_{123}^{\ast}\,\Theta (L)$
induced by the cyclic permutation
$\sigma_{123}$ of the factors of 
$B^3$  transports to an isomorphism
\[ \psi: (\Lambda \times 1)^{\ast}E \simeq \sigma_{123}^{\ast} ((\Lambda \times 1)^{\ast} E )\]
It is convenient to consider this in the fiber above a general element $(x,y,z) \in B^3$. It can
then be displayed as an isomorphism\footnote{The contracted product symbols $\wedge$ will henceforth
mostly be omitted.}
\be
\label{psi} \psi_{x,y,z} : E_{_{x+y,z}}\, E_{x,z}^{-1}\, E_{y,z}^{-1} \simeq E_{_{y+z,x}}\,
E_{y,x}^{-1}\, E_{z,x}^{-1}
\ee Taking into account the maps $E_{z,x}\, E_{z,y} \simeq E_{z,x+y}$ and $E_{x,y}\, E_{x,z} \simeq
E_{x,y+z}$ which expresses the linearity of $E$ in the second variable, one therefore obtains an
isomorphism $\phi$ between the requisite pullbacks of $E$, whose fiber $\phi_{x,y,z}$ above
$(x,y,z)$ is an isomorphism of $A$-torsors
\be
\label{phixyz}
 E_{x+y,z}\,E_{z,x+y}\,E_{y,x}\,E_{x,y} \la E_{y+z,x} \,E_{x,y+z}\,E_{z,y}\,E_{y,z} 
\ee
We impose on the map $\alpha$ (\ref{alph}) the requirement that the following diagram commutes.
\be
\label{phiass}
\xymatrix{
(E_{x+y,z}\,E_{z,x+y})\: (E_{y,x}\,E_{x,y})\ar[rr]^{\phi_{x,y,z}}
\ar[d]_{\alpha_{x+y,z}\,\alpha_{x,y}} && (E_{y+z,x}\,E_{x,y+z}) \:(E_{z,y}\,E_{y,z})
\ar[d]^{\alpha_{x,y+z}\,\alpha_{y,z}}
\\ 
\Lambda(L)_{x+y,z}\:\Lambda(L)_{x,y} \ar[rr] && \Lambda(L)_{x,y+z}\:\Lambda(L)_{x,y}}
\ee
the lower horizontal map being the (obvious) canonical isomorphism. A second constraint on the map
(\ref{alph})  is much simpler to state. It is the requirement that the following commutative
square, in which the horizontal maps are the canonical ones, also commutes:
\be
\label{phicom}
\xymatrix{
\Lambda(L) \ar[r]^{\simeq} \ar[d]_{\alpha} &  s^{\ast}\Lambda(L)\ar[d]_{s^{\ast}\alpha} \\
E \wedge s^{\ast}E \ar[r] & s^{\ast}E \wedge E}
\ee
A final axiom involves both the arrows $\alpha$ (\ref{alph})  and $\beta$ (\ref{bet}). Observe
first of all that the map $\psi$ (\ref{psi}) specializes to
\[\psi_{y,y,x}: E_{_{2y,x}}\, E_{y,x}^{-2} \simeq E_{_{x+y,y}}\, E_{y,y}^{-1}\, E_{x,y}^{-1}\] 
On the other hand, the given $\Sigma$-structure on the first variable of $E$ determines 
an isomorphism \cite{lb:cube} (6.3.2)
\[E_{_{2y,x}} \simeq E^4_{y,x}\]
  so that the map in question determines an isomorphism
\[ E_{_{x+y,y }} \simeq E_{_{x,y}}\, E_{_{y,y}}\, E_{y,x}^2\]
Exchanging the variables $x$ and $y$ this yields 
\[ E_{_{x+y,x }} \simeq E_{_{y,x}}\, E_{_{x,x}} \,E_{x,y}^2\]
These two isomorphisms, and the linearity of $E$ in the second variable therefore determine a
composite map
\[E_{_{x+y,x+y}} \simeq E_{_{x+y, x}} \, E_{_{x+y, y}} \simeq E_{_{x,x}}\, E_{_{y,y}}\,
E^3_{x,y}\, E_{y,x}^3\]
which may also be viewed 
as a map
\be
\xymatrix{
\Lambda\Delta E \ar[r]^<<<<<{\gamma} &E^3 \wedge s^{\ast}E^3}
\ee 
The final condition which we impose on the pair $(E,L)$ is the commutativity of the following
diagram in which the top horizontal arrow is the canonical one
\be
\label{alphabeta1}
\xymatrix{
\Lambda(L^3) \ar[d]_{\Lambda(\beta)} \ar[r]^{\simeq} & \Lambda(L)^3 \ar[d]^{\alpha^3}  \\
\Lambda \Delta E  \ar[r]^<<<<<{\gamma}  &E^3 \wedge s^{\ast}\!E^3}
\ee
\medskip

\noindent We summarize the previous discussion as follows

\begin{definition}

Let $A$ and $B$ be two abelian groups (or sheaves of abelian groups). A $\Gamma_3$-A-torsor pair
above 
$B$ consists of an
$A$-torsor
$L$ on $B$ and an
$A$-torsor $E$ on $B \times$B
satisfying the following conditions:

\medskip

{\it i}) $E$ is endowed with a $\Sigma$-structure with respect to  the first variable and with a
commutative extension structure with respect to  the second. These conditions are compatible,
in the manner previously described  (\ref{compat}).

\smallskip

{\it ii}) There exists a pair of isomorphisms of $A$-torsors 
(\ref{alph}), (\ref{bet})
\[\begin{array}{ccc} 
{\Lambda(L) \stackrel{\alpha}{\la} E \wedge s^{\ast}E} &\hspace{1cm}&\Delta E
\stackrel{\beta}{\la}  L^3
\end{array}
\]

{\it iii}) The diagrams (\ref{phiass}), (\ref{phicom}), and (\ref{alphabeta1}) commute.
\end{definition}

   In order to   obtain a complete 
proof that the category of such quadruples
$(E,L,
\alpha, \beta)$ is described by the truncated complex $t_{\leq 0} \mathrm{RHom}(L\Gamma_3(B),
A[1])$, whose homology groups are the groups  $\mathrm{Ext}^{i}(L\Gamma_3(B), A)$ for $i= 0,1$, one
would have to construct a simple
 realization of the object  $L\Gamma_3(B)$, as was done in  
\cite{lb:cube}
\S 8 for   the object
$L\Gamma_2(B)$. We do not carry this out here, but observe the following fact, which provides, in
view of the exactness of the sequence (\ref{syzg3}),  very strong evidence for this assertion. For
this reason, and for lack of a better name, we will call these quadruples ``$\Gamma_3$-$A$-torsor
pairs on $B$'', or simply ``$\Gamma_3$-torsor pairs''

\begin{proposition}
\label{pr:last}
Let $(E,L,\alpha, \beta)$ be a $\Gamma_3$ torsor pair for which there exists a global section of
$E$, compatible with both the $\Sigma$-structure of $E$ in the first variable and  with the
linearity  in the second one. Then $L$ is endowed with a group law, for which it is  a commutative
extension of
$B$ by
$A$. Furthermore  the pullback of the extension 
\[ \mathrm{0} \la A \la L \la B \la \mathrm{0} \] by the map $B \stackrel{3}{\la} B$ is canonically
split.
\end{proposition}
{\bf Proof:} The given trivialization of $E$ induces by  (\ref{alph}) a trivialization of the
torsor $\Lambda(L)$, and this in turn determines a composition law on $L$ which is compatible with
the projection from $L$ to $B$. The commutativity of diagram (\ref{phiass}) implies that the induced
 trivialization of the source and target  of the lower horizontal map in this diagram 
are compatible, a statement which is equivalent to the associativity of the corresponding
composition law on $L$. Similarly, the commutativity of diagram (\ref{phicom}) implies that the
trivializations of the source and target of the  upper horizontal map in the latter diagram are
compatible, and this in turn is equivalent to the commutativity of the composition law on $L$. This
 shows that $L$ is a commutative extension of $B$ by $A$ (since the
existence of an inverse map on $L$, the existence of a unit element and the compatibility of the
composition law in $L$  with the inclusion of $A$ in $L$ follow automatically). Finally, the
commutativity  of diagram (\ref{alphabeta1}) implies  that the trivialization of the $A$-torsor
$L^3$ determined by the map
$\beta$ (\ref{bet}) is compatible with the group law on $L^3$ induced by the group law of $L$, and
therefore splits the pushout $L^3$ of $L$ by the map $A \stackrel{3}{\la} A$ as a group extension.
This is equivalent to the assertion regarding the pullback of $L$ by the map $B \stackrel{3}{\la}
B$.

\end{document}